\documentclass{article}
\usepackage[latin1]{inputenc}
\usepackage{graphicx}
\usepackage{pgf,tikz}
\usetikzlibrary{arrows}
\usepackage{epstopdf}
\usepackage{dsfont}

\usepackage{amsfonts}
\usepackage{amsmath, latexsym}
\usepackage{amssymb,verbatim}
\usepackage{theorem}
\usepackage{color}
\usepackage{euscript}

\newtheorem{theorem}{Theorem}

\newtheorem{propo}[theorem]{Proposition}

\newtheorem{e-definition}[theorem]{Definition\rm}

\newcommand{\eps}{\varepsilon}

\newcommand{\R}{{\mathbb R}}

\newcommand{\Z}{{\mathbb Z}}

\newcommand{\calL}{{\mathcal{L}}}
\newcommand{\calX}{{\mathcal{X}}}

\newcommand{\dist}{\operatorname{dist}}
\newcommand{\sign}{\operatorname{sign}}

\newcommand{\nd}{\noindent}

\newcommand{\dd}{D}

\newcommand{\vol}{\,\mbox{vol}}

\newcommand{\beq}{\begin{equation}}
\newcommand{\eeq}{\end{equation}}

\title{Interaction energy between vortices of vector fields on Riemannian surfaces}

\author{ 
{\Large Radu Ignat}
\footnote{Institut de Math\'ematiques de Toulouse, Universit\'e
Paul Sabatier, 31062
Toulouse, France. Email: Radu.Ignat@math.univ-toulouse.fr} \and {\Large Robert L. Jerrard}
\footnote{Department of Mathematics, University of Toronto, Toronto, Ontario, Canada. Email: rjerrard@math.toronto.edu}
}

\begin{document}

\maketitle

\begin{abstract}
We study a variational Ginzburg-Landau type model depending on a small parameter $\eps>0$ for (tangent) vector fields on a $2$-dimensional Riemannian surface. As $\eps\to 0$, the vector fields tend to be of unit length and will have singular points of a (non-zero) index, called vortices. Our main result determines the interaction energy between these vortices as a $\Gamma$-limit (at the second order) as $\eps\to 0$.
\end{abstract}

\section{Introduction}
\label{sec:intro}

Let $(S,g)$ be a closed (i.e., compact, connected without boundary) $2$-dimensional Riemannian manifold of genus $\mathfrak g$. We will focus on (tangent) vector fields $$u:S\to TS, \quad \textrm{ i.e., } \,  u(x)\in T_x S \, \textrm{ for every } \, x\in S$$ where $TS=\cup_{x\in S} T_xS$ is the tangent bundle of $S$. It is well known that there are no smooth vector fields $\calX(S)$ (or more generally, of Sobolev regularity $\calX^{1,2}(S)$) of unit length 
$|u|_g=1$ on $S$ (unless $\mathfrak g=1$). In fact, vector fields of unit length have in general singular points with a (non-zero) index. Our aim is to determine the interaction energy between these singular points in a variational model of Ginzburg-Landau type depending on a small parameter $\eps>0$ where the penalty $|u|_g=1$ in $S$ is relaxed. \\

\nd {\it Model}. For vector fields $u:S\to TS$, we define the energy functional
$$
E_\eps(u) 
= \int_S e_\eps(u) \vol_g, \quad e_\eps(u):=
\frac 12 |\dd u|_g^2  +\frac 1{4\eps^2} F(|u|^2_g),
$$
where $|\dd u|_g^2:= |\dd_{\tau_1}u|_g^2 + |\dd_{\tau_2}u|_g^2$ in $S$, $\vol_g$ is the volume $2$-form on $(S,g)$ and $\dd_v$ denotes covariant differentiation (with respect to the Levi-Civita connection) of
$u$ in direction $v$ and $\{ \tau_1,\tau_2\}$ is any local orthonormal basis of $TS$. The potential $F:\R_+\to \R_+$ is a continuous function with $F(1) = 0$ and there exists some $c>0$
such that $F(s^2) \ge c(1- s)^2$ for every $s\geq 0$;  in particular, $1$ is the unique zero of $F$. The parameter $\eps>0$ is small penalizing $|u|_g\neq 1$ in $S$; the goal is to analyse the asymptotic behaviour of $E_\eps$ in the framework of $\Gamma$-convergence (at first and second order) in the limit $\eps\to 0$. This is a ``toy" problem for some physical models arising for thin shells in micromagnetics or nematic liquid crystals (see e.g., \cite{CanSeg,Carbou}).\\

\nd {\it Connection $1$-form}. On an open subset $O\subset S$, a moving frame is a pair of smooth, properly oriented, orthonormal vector fields $\tau_k\in \calX(O)$, $k=1,2$, i.e., $(\tau_k,\tau_l)_g  = \delta_{k\ell}$, $k, l=1,2,$ and $\vol_g(\tau_1,\tau_2) =1$ in $O$,
where $(\cdot, \cdot)_g$ is the scalar product on $TS$. (We will use the same notation $(\cdot, \cdot)_g$ for the inner product associated to $k$-forms, $k=0,1,2$.)
Defining $i:TS\to TS$
such that $i$ is an isometry of $T_xS$ to itself for every $x\in S$ satisfying
\[
i^2w = -w, \qquad
(iw, v)_g = -(w,iv)_g \  = \ \vol_g(w,v), 
\]
then every smooth vector field $\tau \in \calX(O)$ of unit length provides a moving frame $\{\tau_1,\tau_2\} := \{ \tau, i\tau\}$ on $O$. Moreover,  
if $\{\tau_1,\tau_2\}$ is any moving frame in $O$, then $\tau_2=i\tau_1$.
\footnote{In general a moving frame exists only locally on $S$.} Given a moving frame $\{ \tau_1,\tau_2\}$ on an open subset $O\subset S$,
the {\em connection 1-form} $A$ associated to $\{ \tau_1,\tau_2\}$ is defined for every smooth vector field $v\in \calX(O)$:
$$
A(v) := (\dd_v\tau_2,\tau_1)_g = - (\dd_v \tau_1,\tau_2)_g \qquad \textrm{ in }  O.
$$
In particular, $\dd_v \tau_1 = -A(v)\tau_2$ and $\dd_v\tau_2 = A(v)\tau_1$ in $O$.
In complex notation, it yields
for any smooth complex-valued function $\phi$ on $O$:
$$\dd_v( \phi \tau_1 ) =   (d\phi(v) - iA(v) \phi) \tau_1 \qquad \textrm{ in }  O.$$
The definition of $A$ depends on the choice of the moving frame. 
However, the exterior derivative $dA$ of the connection $1$-form is independent of the moving frame, in particular, the following identity holds 
$$
dA = \kappa \ \vol_g,
$$
where $\kappa$ is the Gaussian curvature of $S$
(see \cite{Car94} Proposition 2, Chapter 5.3). We recall
the Gauss-Bonnet theorem that states 
$$
\int_S \kappa \vol_g = 2\pi \chi(S),
$$
where $\chi(S)$ is the Euler characteristic, related to the genus $\mathfrak g$ of $S$ by $\chi(S)=2-2\mathfrak g.$

\medskip

\nd {\it Vortices}. We will identify vortices of a vector field $u$ with small geodesic balls centered at some points around which $u$ has a (non-zero) index. To be more precise, we introduce the Sobolev space $\calX^{1,p}(S)$ of vector fields $u:S\to TS$ such that $|u|_g$ and $|Du|_g$  
belong to $L^p(S)$ (with respect to the volume $2$-form), $p\geq 1$. Given $u\in \calX^{1,p}(S)\cap L^q(S)$ such that $\frac 1p + \frac 1q =1$, $p,q\in [1, \infty]$, we define the $1$-form $j(u)$ 
by \footnote{Note that if $\{\tau_1,\tau_2\}$ is a moving frame on an open set $O\subset S$, then the connection $1$-form $A$ associated to the moving frame is given by $A = - j(\tau_1)$
on $O$. In particular, $d j(u)=-k\vol_g$ in $O$ for every smooth $u\in \calX(O)$ of unit length.} 
$$
j(u)= (\dd u , iu)_g.
$$
In particular, $j(u)$ is a well-defined $1$-form in $L^1(S)$ if $u\in \calX^{1,1}(S)$ with $|u|_g=1$ almost everywhere in $S$; the same is true if $u\in \calX^{1,p}(S)$ for $p\ge \frac 43$. To introduce the notion of index, we assume that $O$ is a simply connected open subset of $S$ and
$u\in \calX^{1,2}(N)$ is a vector field in a neighborhood $N$ of $\partial O$ such that $|u|_g\ge \frac 12$ a.e. in $N$;
then the {\em index} (or winding number) 
of $u$ along $\partial O$ is defined by
$$
\deg (u; \partial O) := \frac 1{2\pi}\left( \int_{\partial O}\frac{ j(u)}{|u|_g^2} + \int_O \kappa \,\vol_g
\right)
$$
(see \cite{Car94} Chapter 6.1). 
In particular, if $u$ is defined in $O$ and has unit length on $\partial O$, then one has $\displaystyle \int_{O} \omega(u)=2\pi \deg (u; \partial O)$ where $\omega(u)$ is the {\it vorticity} associated to the vector field $u$:
\begin{equation}
\omega(u) := d  j(u) + \kappa \vol_g.
\label{omega.def}\end{equation}
Sometimes we can identify the index of $u$ at a point $P\in S$ with the index of $u$ along a curve around $P$.
Note that every smooth vector field $u\in \calX(O)$ (or more generally, $u\in \calX^{1,2}(O)$) of unit length in $O$ has $\deg (u; \partial O)=0$; moreover, a vortex with non-zero index will carry infinite energy $E_\eps$ as $\eps\to 0$.\\

We will prove a $\Gamma$-convergence result (at the second order) of $E_\eps$ as $\eps\to 0$. In particular, at the level of minimizers $u_\eps$ of $E_\eps$, we show that $u_\eps$ converges in $\calX^{1,1}(S)$ (for a subsequence) to a canonical harmonic vector field $u^*$ of unit length that is smooth \footnote{In the case of a surface $(S,g)$ with genus $1$ (i.e., homeomorphic with the flat torus), then $n=0$ and $u^*$ is smooth in $S$. } away from $n=|\chi(S)|$ distinct singular points $a_1, \dots, a_n$, each 
singular point $a_k$ carrying the same index $d_k = \sign \chi(S)$ so that \footnote{In fact, $\deg (u^*; \gamma)=d_k$ for every closed simple curve $\gamma$  around $a_k$ and lying near $a_k$.}  
\beq \label{necessary}
\sum_{k=1}^n d_k=\chi(S).
\eeq
The vorticity $\omega(u^*)$ detects the singular points $\{a_k\}_{k=1}^n$ of $u^*$:  
\beq
\omega(u^*) =  
2\pi \sum_{k=1}^n d_k \delta_{a_k}  \qquad \textrm{ in }  S,
\label{ustar2}
\eeq
where $\delta_{a_k}$ is the Dirac measure (as a $2$-form) at $a_k$. 
The expansion of the minimal energy $E_\eps$ at the second order is given by
$$E_\eps(u_\eps)=n\pi \log \frac 1 \eps+\lim_{r\to 0} \bigg( \int_{S\setminus \cup_{k=1}^n B_r(a_k)}\frac 12 |Du^*|_g^2\, \vol_g + n\pi \log r\bigg) +n\gamma_F+o(1), \textrm{ as } \, \eps\to 0,$$
where $\gamma_F>0$ is a constant depending only on the potential $F$ and $B_r(a_k)$ is the geodesic ball centered at $a_k$ of radius $r$. The second term in the above RHS is called the {\it renormalized energy} between the vortices $a_1, \dots, a_n$ and governs the optimal location of these singular points as in the Euclidian case (see the seminal book \cite{BBH}). In particular, if $S$ is the unit sphere in $\R^3$ endowed with the standard metric $g$, then $n=2$ and $a_1$ and $a_2$ are two diametrically opposed points on $S$. \\

\nd {\it Outline of the note}. The note is divided as follows. Section \ref{sec:canon} is devoted to characterize canonical harmonic vector fields of unit length. In Section \ref{sec:renorm}, we determine the renormalized energy between singular points of canonical harmonic vector fields. The main $\Gamma$-convergence result is stated in the last section. The proofs of these results are part of our forthcoming article \cite{IgnJer}.

\section{Canonical harmonic vector fields of unit length}
\label{sec:canon}

We will say that a canonical harmonic vector field of unit length having the singular points $a_1,\ldots, a_n \in S$ of index
$d_1,\ldots, d_n\in \Z$ for some $n\geq 1$, is a vector field 
$u^* \in \calX^{1,1}(S)$ such that $|u^*|_g=1$ in $S$, \eqref{ustar2} holds and
\beq
d^* j(u^*) =0  \qquad \textrm{ in }  S.
\label{ustar1}
\eeq
Here, $d^*$ is the adjoint of the exterior derivative $d$, i.e., $d^* j(u^*)$ is the unique $0$-form on $S$ such that
$$ \int_S \big(d^* j(u^*), \zeta)_g \vol_g=\int_S \big(j(u^*), d\zeta)_g \vol_g \quad  \textrm{for every smooth $0$-form $\zeta$}.$$ 
If $u^*$ satisfies \eqref{ustar2}, then the 
Gauss-Bonnet theorem combined with \eqref{omega.def} imply that necessarily \eqref{necessary} holds.

We will see that condition \eqref{necessary} is also sufficient. 
Indeed, if \eqref{necessary} holds, we will construct solutions of \eqref{ustar2} and \eqref{ustar1},
as follows:  let $\psi=\psi(a,d)$ be the unique $2$-form on $S$ 
solving:
\begin{equation}
-\Delta \psi = -\kappa\, \vol_g + 2\pi \sum_{k=1}^n d_k\delta_{a_k} \qquad \textrm{ in }  S, \qquad\qquad\int_S \psi = 0,
\label{psi.def}\end{equation}
with the sign convention that $-\Delta = d d^* + d^* d$. The idea is to find $u^*$ such that $j(u^*)-d^*\psi$ is an harmonic $1$-form, i.e., 
$$Harm^1(S)
=\{ \mbox{integrable $1$-forms $\eta$ on $S$} \ : \ d\eta = d^*\eta = 0\mbox{ as distributions}\}.$$
The dimension of the space $Harm^1(S)$ is twice the genus (i.e., $2\mathfrak g$) of $(S, g)$ and we fix an orthonormal basis $\eta_1,\ldots, \eta_{2\mathfrak g}$ of $Harm^1(S)$ such that 
$$
\int_S (\eta_k , \eta_l)_g \ \vol_g= \delta_{kl} \quad \textrm{ for $k,l=1,\ldots, 2\mathfrak g$}.$$
Therefore, it is expected that
\begin{equation}
j(u^*) =   d^*\psi  + \sum_{k=1}^{2\mathfrak g} \Phi_k \eta_k \qquad \textrm{ in }  S
\label{form_jstar}\end{equation}
for some constant vector $\Phi = (\Phi_1,\ldots, \Phi_{2\mathfrak g}) \in \R^{2\mathfrak g}$.
These constants are called {\it flux integrals} as they can be recovered by 
\[
\Phi_k = \int_S (j(u^*), \eta_k)_g \vol_g, \qquad\textrm{ for } k=1,\ldots, 2\mathfrak g.
\]
Note that \eqref{form_jstar} combined with \eqref{psi.def} automatically
yield \eqref{ustar2} and \eqref{ustar1}.
One important point is to characterize
for which values of  $\Phi$  the RHS of \eqref{form_jstar} arises as $j(u^*)$ for
some vector field $u^*$ of unit length in $S$. For that condition, we need to recall the following theorem of Federer-Fleming \cite{FedFle}: there exist $2\mathfrak g$ simple closed geodesics $\gamma_\ell$ on $S$, $\ell=1,\ldots, 2\mathfrak g$, such that for any closed Lipschitz curve $\gamma$ on $S$, one can find integers $c_1\ldots, c_{2\mathfrak g}$
such that 
$$
\textrm{$\gamma$ is homologous to $\displaystyle \sum_{\ell=1}^{2\mathfrak g} c_\ell \gamma_\ell$}$$
i.e., there exists an integrable function $f:S\to \Z$ such that
$$\int_{\gamma}\zeta - \sum_{\ell=1}^{2\mathfrak g} c_\ell\int_{\gamma_\ell}\zeta \ = \ \int_S f \, d\zeta \quad \textrm{ for all
smooth $1$-forms $\zeta$}.$$ 
Having chosen the geodesic curves $\{\gamma_\ell\}_{\ell=1}^{2\mathfrak g}$ and the harmonic $1$-forms $\{\eta_k\}_{k=1}^{2\mathfrak g}$,
we fix the notation
\begin{equation}
\alpha_{\ell k} := \int_{\gamma_\ell}\eta_k, \quad k, \ell=1,\ldots, 2\mathfrak g.
\label{akl.def}\end{equation}

\begin{theorem}\label{P1}
Let $n\geq 1$ and $d = (d_1,\ldots, d_n)\in \Z^n$ satisfy \eqref{necessary}. 
Then for every $a = (a_1,\ldots, a_n)\in S^n$, there exists
\[
\zeta_\ell = \zeta_\ell(a;d) \in  \R/2\pi\Z, \qquad \ell=1,\ldots, 2\mathfrak g
\]
such that if a vector field $u^*\in \calX^{1,1}(S)$ of unit length solves \eqref{ustar2} and \eqref{ustar1},
then 
$j(u^*)$ has the form 
\eqref{form_jstar}  for constants $\Phi_1,\ldots, \Phi_{2\mathfrak g}$ such that
\begin{equation}
\sum_{k=1}^{2\mathfrak g} \alpha_{\ell k} \Phi_k +\zeta_\ell(a,d) \in 2\pi \Z, \qquad \ell=1,\ldots, 2\mathfrak g, 
\label{lattice}\end{equation}
where $(\alpha_{\ell k})$ were defined in \eqref{akl.def}.
Conversely, given any $\Phi_1,\ldots, \Phi_{2 \mathfrak g}$ satisfying 
\eqref{lattice}, there exists a vector field $u^*\in \calX^{1,1}(S)$ of unit length solving \eqref{ustar2} and \eqref{ustar1} and
such that $j(u^*)$ satisfies \eqref{form_jstar}. In addition, the following hold:
\begin{itemize}
\item[1)]  $\zeta_\ell(\cdot ;  d)$ depends continuously on $a\in S^n$ for every $\ell=1,\ldots, 2\mathfrak g$.
More generally, if \footnote{If $\mu$ is a $2$-form (possibly measure-valued) then we write for $p, q\in [1,\infty]$ with 
$\frac 1p + \frac 1q =1$:
\[
\| \mu\|_{W^{-1,p}} := \sup \left\{  \int_{S} f \mu  \ : \ f\in W^{1,q}(S;\R), \,  \| f\|_{W^{1, q}}:=\|f\|_{L^q}+\|df\|_{L^q}\le 1 \right\}.
\]
}
$$
\mu^t := 2\pi\sum_{l=1}^{n_t} d^t_l\delta_{a_l^t}
 \to \mu^0 :=  2\pi\sum_{l=1}^{n_0} d^0_l\delta_{a_l^0}
\qquad\textrm{in $W^{-1,1}$} \qquad\textrm{ as } t\downarrow 0,
$$
$\{d^t_l\}_l$ are integers with \eqref{necessary} and $\sum_{l=1}^{n_t}|d^t_l|$ is uniformly bounded in $t$, 
 then $\zeta_\ell(a^t, d^t)\rightarrow \zeta_\ell(a^0,d^0)$ as $t\downarrow 0$.

\item[2)] any $u^*$ solving \eqref{ustar2} and \eqref{ustar1} belongs to $\calX^{1,p}(S)$ for all $1\leq p<2$,
and is smooth away from $\{a_k\}_{k=1}^n$.

\item[3)] If $u^*, \tilde u^*$ both satisfy \eqref{form_jstar} for the same $(a,d)$ and the same $\{\Phi_k\}_{k=1}^{2\mathfrak g}$, then $\tilde u^* = e^{i\beta}u^*$ for some
$\beta\in \R$.
\end{itemize}
\end{theorem}

\medskip

The constants $\{\zeta_\ell(a;d)\}_{\ell=1}^{2\mathfrak g}$ are determined as follows. For every $\ell=1,\ldots, 2\mathfrak g$, we let 
$\lambda_\ell$ be some smooth simple closed
curve such that $\lambda_\ell$ is homologous to $\gamma_\ell$ (the geodesics
fixed in \eqref{akl.def}) so that $\{a_k\}_{k=1}^n$ is disjoint from
$\lambda_\ell$; for example, $\lambda_\ell$ is either $\gamma_\ell$
or, if $\gamma_\ell$ intersects some $a_k$, a small perturbation thereof.
We now define  $\zeta_\ell(a,d)$ to be the element of  $\R/2\pi\Z$  such that 
\begin{equation}\label{zetak.def}
\zeta_\ell(a,d):=\int_{\lambda_\ell} (d^*\psi + A)  \,   \mod 2\pi, \quad \ell=1,\ldots, 2\mathfrak g,
\end{equation}
where $\psi = \psi(a,d)$ is the $2$-form given by \eqref{psi.def} and $A$ is the connection $1$-form associated to any moving frame defined in
a neighborhood of $\lambda_\ell$. The integral in \eqref{zetak.def} is independent, modulo $2\pi\Z$,
of the choice of moving frame and of the curve
$\lambda_\ell$ homologous to $\gamma_\ell$. 
In examples in which it can be explicitly computed, in general $\zeta_\ell(a,d)\ne 0\mod 2\pi$ for $\ell = 1,\ldots, 2\mathfrak g$.

\section{Renormalized energy}
\label{sec:renorm}

For any $n\geq 1$, we consider $n$ {\bf distinct} points $a = (a_1,\ldots, a_n)\in S^n$. Let $d = (d_1,\ldots, d_n)\in \Z^n$ satisfying \eqref{necessary},  
$\{\zeta_\ell(a;d)\}_{\ell=1}^{2\mathfrak g}$ be given in Theorem \ref{P1} and $\Phi\in \R^{2\mathfrak g}$ be a constant vector inside the set:
$$
\mathcal L(a,d) := \{ \Phi = (\Phi_1,\ldots, \Phi_{2\mathfrak g})\in \R^{2\mathfrak g}\, : \, \sum_{k=1}^{2\mathfrak g} \alpha_{\ell k}\Phi_k +\zeta_\ell(a,d)\in 2\pi \Z, \, \ell=1, \dots, 2\mathfrak g \}.
$$
We define the {\it renormalized energy} between the vortices $a$ of indices $d$ by
$$
W(a,d,\Phi):=\lim_{r\to 0} \bigg(\int_{S\setminus \cup_{k=1}^n B_r(a_k)}\frac12 |Du^*|^2_g\, \vol_g+\pi \log r \sum_{k=1}^n d_k^2\bigg),
$$
where $u^*=u^*(a,d,\Phi)$ is the unique (up to a multiplicative complex number) canonical harmonic vector field given in Theorem \ref{P1} and $B_r(a_k)$ is the geodesic ball centered at $a_k$ of radius $r$. Our arguments show that the above limit indeed exists.
As in the euclidian case (see \cite{BBH}), we can compute the renormalized energy by using the Green's function. For that, let
$G(x,y)$ be the unique function on $S\times S$ such that
$$
-\Delta_x ( G(\cdot, y)  \, \vol_g) = \delta_y \,  - \frac {\vol_g} {\mbox{Vol}_g(S)}\, \, \textrm{ distributionally in } S, \qquad \int_S G(x,y) \vol_g(x) = 0 \quad \mbox{ for every }y\in S.
$$
with $\mbox{Vol}_g(S):= \int_S \vol_g$. 
Then $G$ may be represented in the form (see \cite{Aubin} Chapter 4.2):
$$
G(x,y) = G_0(x,y) + H(x,y),\qquad\mbox{ with }H\in C^1(S\times S),
$$
where $G_0$ is smooth away from the diagonal, 
 with
$$
G_0(x,y) = -\frac 1{2\pi}\log(dist(x,y)) \, \,  \mbox{ if the geodesic distance } dist(x,y) <\frac 12 (\mbox{injectivity radius of $S$}).
$$
The $2$-form $\psi=\psi(a,d)$ defined at \eqref{psi.def} can be written as: 
$$
\psi= 2\pi \sum_{k=1}^n d_k G(\cdot, a_k) \vol_g \ +  \psi_0  \vol_g \qquad\mbox{ in } S,
$$
where $\psi_0\in C^\infty(S)$ has zero average on $S$ and solves
\begin{equation}
-\Delta \psi_0  = -\kappa + \bar \kappa , \qquad
\mbox{ for }
\bar \kappa = \frac 1{\mbox{Vol}(S)} \int_S \kappa \vol_g  = \frac {2\pi \chi(S)}{\mbox{Vol}(S)} .
\label{psi0.def}\end{equation}
In other words, the $2$-form $x\mapsto \psi(x) + d_k \log \dist(x, a_k) \vol_g$ is $C^1$ in a neighborhood of $a_k$ for every $1\leq k\leq n$. We have the following expression of the renormalized energy:

\medskip

\begin{propo}
Given $n\geq 1$ distinct points $a_1,\ldots, a_n \in S$, integers $d_1,\ldots, d_n$ with \eqref{necessary}
and $\Phi \in \mathcal L(a,d)$, then
\beq
W(a,d,\Phi) =  4\pi^2 \sum_{l\ne k} d_l d_k G(a_l,a_k) 
+2\pi \sum_{k=1}^n  \left[ \pi d_k^2 H(a_k,a_k) +d_k \psi_0(a_k) \right]  
+ \frac 12|\Phi|^2 + \int_S \frac{|d\psi_0|^2}2 \vol_g \ ,
\label{W.formula}\eeq
where $\psi_0$ is defined in \eqref{psi0.def}. 
\label{prop.W}\end{propo}

\medskip

In the case of the unit sphere $S$ in $\R^3$ endowed with the standard metric (in particular, $\psi_0$ vanishes in $S$), 
if $n=2$ and $d_1=d_2=1$, then the second term in the RHS of \eqref{W.formula} is independent of $a_k$ (as $x\mapsto H(x,x)$ is constant, 
see \cite{steiner}); moreover, $\Phi=0$ and so, minimizing $W$ is equivalent by minimizing the Green's function $G(a_1, a_2)$ over the set of 
pairs $(a_1, a_2)$ in $S\times S$, namely, the minimizing pairs are diametrically opposed. 

\section{$\Gamma$-convergence}

Given the potential $F$ in Section \ref{sec:intro}, we compute the energy $E_\eps$ of the radial profile of a vortex of index $1$ inside a geodesic ball of radius $R>0$:
$$
I_F(R,\eps) := \inf \left \{ \pi \int_0^{R}\left[ f'(r)^2 + \frac {f(r)^2}{r^2} + \frac{1}{2\eps^2}F(f(r)^2)\right] r dr \ 
: f(0)=0, f(R) = 1 \right \}.
$$
Then $I_F(R,\eps) = I_F(\lambda R, \lambda \eps) = I_F(1,\frac{\eps}{R})=:I_F(\frac{\eps}{R})$ for every
$\lambda>0$, and the following limit exists (see \cite{BBH}): 
$$
\gamma_F:=\lim_{t\to 0}( I_F(t)  + \pi \log t).
$$

We state our main result: \\

\begin{theorem}\label{intrinsic.gammalim}
The following $\Gamma$-convergence result holds.
\begin{itemize}
\item[1)] (Compactness) Let $(u_\eps)_{\eps\downarrow 0}$ be a family of vector fields on $S$ satisfying
$
E_\eps(u_\eps)  \le N \pi |\log \eps| + C
$
for some integer $N\geq 0$ and a constant $C>0$. Denoting by
$$\Phi(u_\eps) := \left(\int_S  (j(u_\eps) ,  \eta_1)_g \vol_g ,\ldots, \int_S  (j(u_\eps), \eta_{2\mathfrak g})_g \vol_g\right)\in \R^{2\mathfrak g},$$
then there exists a sequence $\eps \downarrow 0$ such that 
\begin{equation}
\omega(u_\eps)  \longrightarrow 2\pi \sum_{k=1}^{n} d_k\delta_{a_k} \quad \textrm{in } \, W^{-1,1}, \quad \Phi(u_\eps)\rightarrow \Phi \, \,  \textrm{ as } \, \eps \to 0,
\label{convergence}\end{equation}
where $\{a_k\}_{k=1}^n$ are distinct points in $S$ and $\{d_k\}_{k=1}^n$ are nonzero integers satisfying \eqref{necessary} and $\sum_{k=1}^n |d_k|\leq N$ and $\Phi\in \calL(a,d)$. Moreover, if $\sum_{k=1}^n |d_k|= N$, then
$n = N$ and 
$|d_k|=1$ for every $k=1, \dots, n$ (in particular, $n=\chi(S)$ modulo $2$).

\item[2)] ($\Gamma$-liminf inequality) Assume that the vector fields $u_\eps\in \calX^{1,2}(S)$ satisfy \eqref{convergence} for $n$ distinct points $\{a_k\}_{k=1}^n\in S^n$ and $|d_k|=1$, $k=1, \dots n$ that satisfy \eqref{necessary} and $\Phi\in \calL(a,d)$. Then
$$\liminf_{\eps\to 0}
\left[
E_\eps(u_\eps)- n \pi |\log \eps|) 
\right]  \ \ge \  
W(a,d,\Phi)+n \gamma_F.$$ 

\item[3)] ($\Gamma$-limsup inequality) For every $n$ distinct points $a_1,\ldots, a_n\in S$ and $d_1,\ldots , d_n\in \{\pm 1\}$ satisfying
\eqref{necessary} and every $\Phi\in \calL(a,d)$
there exists a sequence of vector fields $u_\eps$ on $S$ such that \eqref{convergence} holds and
$$
E_\eps(u_\eps)- n \pi |\log \eps| \longrightarrow  W(a,d, \Phi)+n \gamma_F \quad \textrm{ as } \, \eps \to 0.
$$
\end{itemize}
\end{theorem}

This theorem is the generalization of the $\Gamma$-convergence result for $E_\eps$ in the euclidian case (see {\cite{CoJe,JeSo,SandSerfbook,AliPon}) and it is based on topological methods for energy concentration (vortex ball construction, vorticity estimates etc.) as introduced in \cite{Je,Sa}.

\section*{Acknowledgements}
R.I. acknowledges partial support by the ANR project ANR-14-CE25-0009-01. 
The research of R.J. was partially supported by the National Science and Engineering Research Council of Canada under operating grant 261955.

\end{document}